\documentclass[10pt,english,reqno]{amsart}
\usepackage{amsmath,amssymb,amsthm,amsfonts,graphicx,color}
\usepackage{amssymb,geometry}

\usepackage{color, graphics}
\usepackage{dsfont}
\usepackage{graphicx}
\usepackage[T1]{fontenc}
\usepackage[latin1]{inputenc}
\usepackage[english]{babel}
\usepackage{geometry}
\usepackage{cases}
\usepackage[colorlinks=true]{hyperref}

\usepackage[colorlinks=true]{hyperref}
\hypersetup{linkcolor=black,citecolor=black,filecolor=black,urlcolor=black}

\newcommand{\h}{{\tt h}}

\newcommand{\pe}{{\tt p}}

\newcommand{\x}{{\tt x}}

\newcommand{\R}{\mathbb{R}}

\newcommand{\diver}{\text{div}}

\newcommand{\eps}{\varepsilon}

\newtheorem{definition}{Definition}[section]
\newtheorem{theorem}{Theorem}[section]

\newtheorem{remark}{Remark}[section]
\newcommand{\bremark}{\begin{remark} \em}
\newcommand{\eremark}{\end{remark} }

\numberwithin{equation}{section}

\title{Lawson cones and the Allen-Cahn equation}

\author{Oscar Agudelo}
\address{University of West Bohemia in Pilsen-NTIS, Univerzitn\'{i} 22, Czech Republic.}
\email {oiagudel@ntis.zcu.cz}

\author{Matteo Rizzi}
\address{Mathematisches Institut, Justus Liebig Universit\"{a}t, Arndtstrasse 2, 35392, Giessen, Germany.}
\email{mrizzi1988@gmail.com}

\begin{document}
\maketitle

\begin{abstract}
In this paper we discuss nondegeneracy and stability properties of some special minimal hypersurfaces which are asymptotic to a given Lawson cone $C_{m,n}$, for $m,\,n\ge 2$. Then we use such hypersurfaces to construct solutions to the Allen-Cahn equation $-\Delta u=u-u^3$ in $\R^{N+1}$, $N+1\ge 8$, whose zero level set has exactly $k\ge 2$ connected components and with infinite Morse index.
\end{abstract}

\section{Introduction}

In this work we discussed part of the recent developments related to entire and bounded solutions to 
the {\it Allen-Cahn equation}
\begin{equation}\label{allen-cahn}
\Delta u-F'(u)=0\quad\hbox{in}\quad  \R^{N+1},
\end{equation}
$F(s)=\frac{1}{4}(1-s^2)^2$ for $s\in \R$ is the best example of a {\it symmetric double well potential}.

Equation \eqref{allen-cahn} appears in the modeling of phase transitions phenonema (see \cite{Allen-Cahn}). Roughly speaking, the functions $u_+=1$ and $u_-=-1$ represent two different phases of a material placed in a container. When the two phases co-exist (not mixing involved), they are separated by an interface. The equilibrium configuration is the smooth function $u$, which describes the phase change, and is expected to be a solution of \eqref{allen-cahn}. Observe that $u_{\pm}$ are trivial solutions of \eqref{allen-cahn} corresponding to stable phases of the said material.

Theoretical developments related to the equation \eqref{allen-cahn} are known to be strongly connected with the {\it Theory of Minimal Surfaces}. This statement can be briefly explained as follows. First, we remark that in dimension one, the function 
\begin{equation}\label{heteroclinic}
w(z)=\tanh\left(\frac{z}{\sqrt{2}}\right), \quad z\in \R.
\end{equation}
is the unique solution (up tu translation of the argument) to the problem
\begin{equation}\label{eqn:1D_AC}
w'' +w(1-w^2)=0, \quad \hbox{in } \R,\quad  w(\pm \infty)=\pm1, \quad w'>0.
\end{equation}

It is verified that
$$
\small{w(z)=\left\{
\begin{array}{ccc}
1 -2 e^{-\sqrt{2}\,z} + O\left(e^{-2\sqrt{2}\,z}\right), \quad z>>0\\
-1 +2 e^{-\sqrt{2}\,z} + O\left(e^{-2\sqrt{2}\,z}\right), \quad z<<0.\\
\end{array}
\right.}
$$

Next, the equation \eqref{allen-cahn} appears as the limit of its singularly perturbed version, namely
\begin{equation}\label{ScaledAllen-CahnEqn}
\epsilon^2 \Delta v-F'(v)=0 \quad \hbox{in } \Omega, \qquad\
\frac{\partial v}{\partial \nu}=0\quad \hbox{on }\partial \Omega, \qquad  \epsilon\to 0,
\end{equation}
where $\Omega\subset \R^{N+1}$ is a smooth bounded domain. As discussed in \cite{Mo}, given a solution $v$ of \eqref{ScaledAllen-CahnEqn}, the rescaling $u_{\epsilon}(x)=v(\epsilon x)$ yields that 
$$
\Delta u_{
\epsilon
} - F'(u_{\epsilon})=0 \quad \hbox{in} \quad \Omega_{\epsilon}:=\epsilon^{-1}\Omega.
$$

Since, $\epsilon^{-1}\Omega \to \R^N$, from standard Schauder estimates, in the limit we obtain a smooth bounded function $u
$ that solves \eqref{allen-cahn}.

\medskip
On the other hand, regarding the solutions of the equation \eqref{ScaledAllen-CahnEqn}, they are exactly the critical points in $H^1(\Omega)$ of the energy
$$
J_{\epsilon}(v)=\int_{\Omega} \frac{\epsilon}{2}\,|\nabla v|^2 +
\frac{1}{\epsilon}\,F(v).
$$

Observe that $v=\pm1$ are global minimizers of $J_{\epsilon}$. The capture of configurations with coexisting phases is based on the intuition that such configurations $v_{
\epsilon}$ should resemble an $\epsilon-$regularization of the function 
$$ v^* = \mathds{1}_{\Lambda} - \mathds{1}_{\Omega \setminus\Lambda},$$
where $\Lambda
\subset \Omega$. The characteristic functions $\mathds{1}_{\Lambda}$ and $\mathds{1}_{\Omega \setminus \Lambda}$ of $\Lambda$ and $\Omega\setminus \Lambda$ respectively  might be thought as  constant phases of a material. Observe also that $v^*$ minimizes the second term in $J_{\epsilon}$, but $v^*\notin H^1(\Omega)$ (see \cite{Mo}).

Let $M$ denote the boundary of $\Lambda$ relatively to $\overline{\Omega}$. The ideas in \cite{Mo} are based on the intuition that when $M$ is a smooth orientable hypersurface, the transition from the phase $\mathds{1}_{\Lambda}$ to the phase $\mathds{1}_{\Omega\setminus\Lambda}$ takes place along the normal direction of $M$ with a approximate one dimensional profile, i.e.,  
$$ v_{\epsilon}(x)\approx
w\left(\frac{z}{\epsilon}\right),
$$
where $z$ corresponds to the normal direction to $M$ and $w$ is described in \eqref{heteroclinic} and \eqref{eqn:1D_AC}. The gradient term in $J_{\epsilon}(v_{\epsilon})$ suggests that $M$ should be critical for the area. To make this idea slightly more precise, we take $\epsilon\to 0^+$, so that
$$
J_{\epsilon}(v_{\epsilon}) \approx {\rm Area}(M)\int_{\R}\left(\frac{1}{2}|w'|^2 + F(w)\right)dz.
$$

This intuition lead to
Modica's result in \cite{Mo}, which states that for a family of local minimizers $\{v_{\epsilon}\}_{\epsilon}$ of $J_{\epsilon}$ with bounded energy, 
$$
v_{\epsilon}\to v^* \quad \hbox{in} \quad L^1_{loc},
$$
where $M=\partial \Lambda$ has minimal
perimeter. We remark that in case $M$ is a smooth hypersurface, this means that $M$ is a minimal hypersurface.  Observe also that $M$ is the approximate nodal set of $v_{\epsilon}$ for $\epsilon>0$ small.

Kohn and Sternberg in \cite{KS} constructed a family $\{v_\eps\}_\eps\in(0,\eps_0)$ of solutions of \eqref{ScaledAllen-CahnEqn} on a bounded domain $\Omega\subset\R^2$ whose zero level set is approaching some given segment as $\eps\to 0$. A similar result on compact manifolds was proved by Pacard and Ritor\'{e} \cite{PR} on more general manifolds.

Modica's developments and many other related developments during the 70th's, lead De Giorgi to formulate its remarkable conjecture in 1978.\\

{\bf De Giorgi conjecture (1978):} Let $u\in L^{\infty}(\R^N)$ be an  entire solution of the equation. Assume also that $u$ is \emph{monotone in one direction} and that $2\leq N+1\leq 8$. Then, after rotation and translation $u(x)=w(x_1)$,  where is determined by \eqref{eqn:1D_AC}.\\

Some of the key developments in proving or disproving the De Giorgi's conjecture are the following. This is motivated by the fact that entire minimal graphs over $\R^{N+1}$ are affine if $N\le 7$ ({\it Bernstein conjecture}
). The conjecture is known to be true in dimension $N+1=2,\,3$ (see \cite{GG,AC}). In dimension $4\le N+1\le 8$ only partial results are available. For instance Savin \cite{Sa} proved that in those dimensions the conjecture is true under the further assumption that $u(x_1,x')\to\pm 1$ as $x_1\to\pm\infty$, for any $x'\in\R^N$. Moreover, it is known that the conjecture is sharp about the dimension. In fact in \cite{DKW-degiorgi} the authors constructed a family of monotone solutions in $\R^9$ which are not one-dimensional.

At the core of the De Giorgi's conjecture is the hypothesis of monotonicity of the solutions, which implies their stability. To be more precise, let $u\in L^{\infty}(\R^N)$ be a solution of \eqref{allen-cahn} and consider the linearized operator of \eqref{allen-cahn} around $u$, namely 
$$
L\psi:=\Delta \psi + (1-3u^2)\psi\quad \hbox{in }\R^N, \qquad \psi \in L^{\infty}(\R^N).
$$

The stability of $u$ is equivalent to saying that the operator $L$ satisfies the maximum principle or, in turn it is equivalent to say that the associated quadratic form 
$$
B(\psi,\psi):=\int_{\R^N} |\nabla \psi|^2 - (1-3u^2)\psi^2.
$$
is positive definite in the space $C_{c}^{\infty}(\R^N)$.

The {\it Morse Index of u}, $m(u)$, is the  maximal dimension of a vector space of smooth functions $\psi$ with \emph{compact support} for which
$$
B(\psi,\psi)<0.
$$

A solution $u\in L^{\infty}(\R^N)$ of \eqref{allen-cahn} is unstable if $m(u)>0$. 

Regarding unstable solutions the following are some of the recent and relevant examples and results. Dang, Fife, Peletier in \cite{DFP} exhibited saddle solutions in dimension two. Their result was later generalized to higher dimensions by Cabre and Terra, first in 2010 by showing unstable saddle solutions and then in 2011 by showing stable saddle solutions (see \cite{CT2009,CT2010}).

In 2011, del Pino, Kowalczyk, Wei in \cite{DKW} constructed solutions to \eqref{allen-cahn} in three dimensions. The nodal set of such solutions is close to a large dilation of a complete embedded
non-degenerate minimal surface with finite total curvature. These solutions also have finite Morse Index. As a follow up for the three dimensional setting, in 2012 Agudelo, del Pino and Wei constructed two types of solutions to \eqref{allen-cahn}, one with Morse Index one and disconnected nodal set asymptotically characterized by Liouville's equation. The other one with Morse index going to infinity and nodal set characterized by nested catenoids with logarithmically diverging ends (see \cite{ADW2015}). These solutions correspond to limiting situations of a single one-parameter family of solutions (see \cite{GLW}).

Concerned with the higher dimensional case, in \cite{ADW2016}, the same authors constructed solutions  to \eqref{allen-cahn} whose nodal set resembles a $N-$catenoid  with a logarithmic correction at infinity. In \cite{Cabre2011}, Cabr\'e showed stable solutions of \eqref{allen-cahn} in dimension $N+1=2m$ with $m\leq 7$. The nodal set of these solutions is the {\it Simon's cone} $C_m$. In 2013, Pacard and Wei construct stable solutions to \eqref{allen-cahn} in dimension $N+1$ with $N\geq 7$. The nodal set of these solutions is close to a foliated {\it Simon's cone} like  smooth surface.

Summarizing, the classification and the construction of solutions to the Allen-Cahn equation \eqref{allen-cahn} is strongly related to the theory of minimal hypersurfaces, that is hypersurfaces with zero mean curvature. In this work, we present a recent contribution that is a follow up of the results mentioned above. To be more precise, we discuss stability and nondegeneracy results for some special minimal hypersurfaces, which are asymptotic to a specific cone (see Section \ref{sec-geom}) and we also present some recent results concerned with existence and stability of solutions to \eqref{allen-cahn} (see Sect 3).\\

The plan of the paper is the following: in Section \ref{sec-geom} we deal with stability and non-degeneracy properties of some special minimal hypersurfaces, in Section \ref{sec-Allen-Cahn} we use these hypersurfaces to construct entire solutions to the Allen-Cahn equation with infinite Morse index
whose zero level set is the union of $k\ge 2$ normal graphs over one of these hypersurfaces.

\section{minimal hypersurfaces: stability and nondegeneracy}\label{sec-geom}

In this part we discuss further some of the underlying geometric setting related to hysurfaces with minimal perimeter. 

We define the perimeter of a measurable subset $E\subset\R^{N+1}$ in an open set $\Omega\subset\R^{N+1}$ as
\begin{equation}\label{Perimeter}
{\rm Per}(E,\Omega):=\sup\left\{\int_E \diver Y \,d\xi:\, Y\in C^\infty_c(\Omega,\R^{N+1}),\,|Y|\le 1\right\}.
\end{equation}
If $E$ has smooth boundary $\partial E$, it follows from the Divergence Theorem that ${\rm Per}(E,\Omega)$ coincides with the $N$-dimensional {\it Hausdorff measure} of $\partial E\cap\Omega$. However, the definition in \eqref{Perimeter} allows us to treat the case of sets $E$ with non-smooth boundary $\partial E$.\\

At this point we give precise context of the conclusion of Modica's results stated in \cite{Mo}, that is, the set $\Lambda:=\{x\in \Omega:\,v^*(x)=1\}$ solves the minimisation problem
$${\rm Per}(\Lambda,\Omega)=\min\left\{{\rm Per}(F,\Omega):|F|=\frac{m+|\Omega|}{2}\right\},$$
where $$m=\int_\Omega v^*dx\in(-|\Omega|,|\Omega|).$$
Let $\Sigma \subset \R^{N+1}$ be a hypersurface with singular set ${\rm sing}(\Sigma)$. Assume that ${\rm sing}(\Sigma)$ has zero $N$-dimensional Hausdorff measure and that $\Sigma$ is orientable with $\nu_\Sigma:\Sigma\setminus{\rm sing(\Sigma)} \to {S}^{N}$ a continuous choice of its unit normal vector field. We set
$${\rm Area}(\Sigma):=\int_\Sigma d\sigma,$$
where $d\sigma$ is the $N$-dimensional Hausdorff measure. We note that ${\rm Area}(\Sigma)\in[0,\infty]$ and it coincides with the usual notion of area when $\Sigma$ is smooth. Taking a function $\phi\in C^\infty_c(\Sigma\setminus{\rm sing}(\Sigma))$ and considering, for $\eps$ small enough, the normal variation 
$$\Sigma_\eps:=\{y+\eps\phi(y)\nu_\Sigma(y):\,y\in{\rm supp}(\phi)\}$$
and differentiating in $\eps$, we have
\begin{equation}\label{eqn:der_area_sigma_eps}
\frac{d}{d\eps} {\rm Area}(\Sigma_\eps)\bigg|_{\eps=0}=-\int_\Sigma H_\Sigma \phi d\sigma,
\end{equation}
where $H_{\Sigma}:= \kappa_1 + \cdots + \kappa_N$ is the mean curvature of $\Sigma$ and $\kappa_1,\ldots,\kappa_N: \Sigma \setminus {\rm sing}(\Sigma)\to \R$ correspond to the principal curvatures of $\Sigma$. The hypersurface $\Sigma$ is said to be \textit{minimal} if $H_\Sigma\equiv 0$ in $\Sigma\backslash{\rm sing}(\Sigma)$ or equivalently, since the function $\phi$ is arbitrary, if $\Sigma$ is a critical point of the Area functional with respect to compactly supported normal variations. Observe that, according to our definition, minimal hypersurfaces are only \textit{critical points} of the Area functional, not necessarily minimisers. The true minimisers of the Area functionals are called \textit{Area minimising} hypersurfaces. More precisely, taking an open set $E\subset\R^{N+1}$, a bounded open set $\Omega\subset\R^{N+1}$ and assuming, without loss of generality, that $0\in\partial E$, we introduce the following definition.
\begin{definition}
\label{def_area-minimising}
We say that $\Sigma:=\partial E$ is an area-minimising (or minimising) hypersurface if for any $\rho>0$ and for any smooth set $F\subset\R^{N+1}$ such that $F\backslash B_\rho(0)=E\backslash B_\rho(0)$,
$$
{\rm Per}(E,B_{2\rho}(0))\le {\rm Per}(F,B_{2\rho}(0)).
$$
\end{definition}
An area minimising hypersurface is minimal, but the opposite implication is not true. Therefore, in order to give a more accurate description of minimal hypersurfaces, it is crucial to consider the second variation of the ${\rm Area}$ functional. More precisely, in the above notations for the normal variation of $\Sigma$, we introduce the second variation of ${\rm Area}(\Sigma)$ as follows
\begin{equation}
\label{second-variation-Area}
\mathcal{Q}_\Sigma(\phi):=\frac{d^2}{d\eps^2}{\rm Area}(\Sigma_\eps)\bigg|_{\eps=0}=\int_\Sigma \big(|\nabla_\Sigma \phi|^2-|A_\Sigma|^2\phi^2) \, d\sigma.
\end{equation}
We refer the reader to \cite{Nunes} for the detailed calculations.
\begin{definition}
\label{def-Sigma-stable}
Let $\Sigma\subset \R^{N+1}$ be a minimal hypersurface of codimension one with ${\rm sing}(\Sigma)$ having measure zero and such that $\Sigma\setminus {\rm sing}(\Sigma)$ is an orientable hypersurface. Assume that $\Sigma$ is minimal. We say that
\begin{enumerate}
\item $\Sigma$ is stable if for any $\phi\in C^\infty_c(\Sigma\setminus{\rm sing(\Sigma)})$, $\mathcal{Q}_\Sigma(\phi)\ge 0$;  
\item $\Sigma$ is strictly stable if $$\inf\left\{\mathcal{Q}_\Sigma(\phi):\,\phi\in C^\infty_c(\Sigma\setminus{\rm sing}(\Sigma)),\,\||A_\Sigma|\phi\|_{L^2(\Sigma)}=1\right\}>0.$$
\end{enumerate}
\end{definition}
Integrating by parts it is possible to see that
$$\mathcal{Q}_\Sigma(\phi)=\int_\Sigma -J_\Sigma \phi\,\phi\, d\sigma\qquad\forall\,\phi\in C^\infty_c(\Sigma\setminus {\rm sing}(\Sigma)),$$
where $J_\Sigma:=\Delta_\Sigma+|A_\Sigma|^2$ is known as the \textit{Jacobi operator} of $\Sigma$, $\Delta_\Sigma$ is the Laplace-Beltrami operator of $\Sigma$ and $|A_\Sigma|^2:=\kappa_1^2+\dots+\kappa_N^2$ is the squared norm of the second fundamental form. A distributional solution $\phi$ to the Jacobi equation 
\begin{equation}
\label{Jacobi-eq}
    J_\Sigma\phi=0
\end{equation}
is known as a \textit{Jacobi field} of $\Sigma$. We will be mostly interested in bounded Jacobi fields, but singular or unbounded Jacobi fields will play also a role in our discussion. 

Next, let us introduce the Lawson cone
\begin{equation}\label{def:Lawson_cone}C_{m,n}:=\{(x,y)\in\R^m\times\R^n:\,(n-1)|x|^2=(m-1)|y|^2\},\qquad m,n\ge 2, \qquad N+1=m+n.
\end{equation}

This cone is minimal, since its mean curvature is zero everywhere, except the origin. Also, $C_{m,m}$ is known as the Simon's cone.

We are interested in some special hypersurfaces $\Sigma$, which are asymptotic to these cones. In this regard and in high dimension we have the following result.

\begin{theorem}\label{th_Al}{\cite{AR2}}
Let $m,n\ge 2$, $m+n\ge 8$. Then there exist exactly two smooth minimal hypersurfaces $\Sigma^\pm_{m,n}\subset E^\pm_{m,n}$ satisfying that
\begin{enumerate} 
\item $\Sigma^\pm_{m,n}$ are asymptotic to $C_{m,n}$ at infinity.
\item ${\rm dist}(\Sigma^{\pm}_{m,n},\{0\})=1$;
\item  $\Sigma^\pm_{m,n}$ are $O(m)\times O(n)$-invariant.
\end{enumerate}
Moreover, the Jacobi field $y\cdotp\nu_{\Sigma^\pm_{m,n}}(y)$ never vanishes. 
\end{theorem}
For these hypersurfaces we proved the following stability result (see \cite{AR3}).
\begin{theorem}[\cite{AR3}]
\label{th-strict-stability}
The hypersurfaces constructed in Theorem \ref{th_Al} are strictly stable.
\end{theorem}
The stability of the hypersurfaces constructed in Theorem \ref{th_Al} directly follows from point $(3)$ of Theroem \ref{th_Al}. The nontrivial point of Theorem \ref{th-strict-stability} is that the stability is actually strict. This follows from an injectivity result for the the Jacobi operator $J_\Sigma$ of $\Sigma$ in some suitable function spaces (see Sections $3$ and $4$ of \cite{AR3} for the details). In any case, strict stability does not rule out the existence of bounded and even decaying Jacobi fields, such as the Jacobi fields $\{\nu_\Sigma(y)\cdotp e_i\}_{1\le i\le N+1}$ and the Jacobi field $y\cdotp\nu_\Sigma(y)$ coming from dilations. Strict stability implies that the Jacobi fields do not have fast decay at infinity.

Next, we consider hypersurfaces which are asymptotic to a fixed Lawson cone in lower dimension. 
\begin{theorem}[\cite{ABPRS, M}]
\label{th-Sigma-low-dim}
Let $m,n\ge 2$, $m+n\le 7$. Then there exists a unique complete, embedded, $O(m)\times O(n)$ invariant minimal hypersurface $\Sigma_{m,n}$ such that 
\begin{enumerate}
\item \label{Sigma-as} $\Sigma_{m,n}$ is asymptotic to the cone $C_{m,n}$ at infinity,
\item \label{Sigma-osc} $\Sigma_{m,n}$ intersects $C_{m,n}$ infinitely many times,
\item \label{Sigma-ort} $\Sigma_{m,n}$ meets $\R^m\times\{0\}$ orthogonally,
\item \label{Sigma-dist} $dist(\Sigma_{m,n}, C_{m,n})=1$.
\end{enumerate}
\end{theorem}
For such hypersurfaces it is relevant to compute the Morse index, defined as
\begin{equation}
Morse(\Sigma):=\sup\{\dim(X):\,X\text{ subspace of }C^\infty_c(\Sigma\setminus{\rm sing}(\Sigma)):\,\mathcal{Q}_\Sigma(\phi)<0,\,\forall\,\phi\in X\backslash\{0\}\}.
\end{equation}

It follows from the definition that $Morse(\Sigma)=0$ if and only if $\Sigma$ is stable. For our hypersurfaces we have the following result.
\begin{theorem}[\cite{AR3}]
\label{th-Morse-infinite}
The hypersurfaces constructed in Theorem \ref{th-Sigma-low-dim} have infinite Morse index.
\end{theorem}
The proof of Theorem \ref{th-Morse-infinite} given in Section $4$ of \cite{AR3} is explicit. In fact it is possible construct a infinite-dimensional subspace of compactly supported functions which make the quadratic form $\mathcal{Q}_\Sigma$ negative. This is possible thanks to the fact that $3\le N\le 7$.

Minimal hypersurfaces are invariant under dilations, translations and rotations. Such geometric transformations give rise to at most $\frac{N(N+1)}{2}+N+2$ linearly independent Jacobi fields, which are known as \textit{geometric Jacobi fields}.
\begin{definition}
\begin{itemize}
\item A minimal hypersurface $\Sigma\subset\R^{N+1}$ is said to be nondegenerate if all its bounded Jacobi fields are geometric.
\item Let $S\subset O(N+1)$ be a subgroup. We say that $\Sigma$ is $S$-nondegenerate if all $S$-invariant Jacobi fields are geometric.
\end{itemize}
\end{definition}
A similar definition was given, for instance, in \cite{DKW} for the Costa-Hofmann-Meeks surfaces. In particular, for our hypersurfaces, we have the following nondegeneracy result.
\begin{theorem}
\label{th-nondegeneracy}
\begin{enumerate}
\item The hypersurfaces constructed in Theorem \ref{th_Al} are nondegenerate.
 \item The hypersurface $\Sigma_{m,n}$ constructed in Theorem \ref{th-Sigma-low-dim} is $O(m)\times O(n)$-nondegenerate, 
 for any $m,n\ge 2$, $m+n\le 7$.
\end{enumerate}
\end{theorem}
The proof of this nondegeneracy result is based on a Fourier expansion of the Jacobi fields of our hypersurfaces and either the strict stability in case $N+1\ge 8$ or the symmetry properties for $4\le N+1\le 7$, which allow us to reduce the Jacobi equation $J_\Sigma\phi=0$ to a second order ODE.

\section{The Allen-Cahn equation}\label{sec-Allen-Cahn}

In this section we turn again up to equation \eqref{allen-cahn} and we discuss the connection with the Lawson cone $C_{m,n}$ introduced in \eqref{def:Lawson_cone}. We observe that $C_{m,n}$ is invariant under the action of the group $O(m)\times O(n)$ and that $\R^{N+1}\backslash C_{m,n}$ has {\it two connected components}, namely, 
$$
\begin{aligned}
E^+ &:=\bigg\{(x,y)\in\R^m\times\R^n:\,|x|^2<\frac{m-1}{n-1}|y|^2\bigg\} \quad {\hbox{inside of } C_{m,n};} \\ 
E^-& :=\bigg\{(x,y)\in\R^m\times\R^n:\,|x|^2>\frac{m-1}{n-1}|y|^2\bigg\} \quad {\hbox{outside of } C_{m,n}.}
\end{aligned}
$$

Assume  $N=m+n\geq 8$ and let $
\Sigma_{m,n}^+\subset E^+$ and $\Sigma_{m,n}^-\subset E^-$ be the two smooth stable minimal surfaces discussed in Theorem \ref{th_Al}. In what follows $\Sigma$ denotes either $\Sigma^+_{m,n}$ or $\Sigma^-_{m,n}$ and $\nu:\Sigma \to S^{N}$ is a fixed choice of the unit normal vector.

Next, for $\eps>0$ small enough, set $\Sigma_{\eps}:=\eps^{-1}\Sigma$ and consider the associated {\it Fermi coordinates}
$$
\x=\pe + z\nu(\eps \pe)
$$
for $\pe\in \Sigma_{\eps}$ and $|z|<\frac{\delta}{\eps}$. The first existence result that we present is the following.
  
\begin{theorem}[\cite{AKR}]
\label{Theorem 1} Let $m,\,n\ge 3$, $n+m=N\ge 8$ and let $\Sigma$ be one of the hypersurfaces constructed in Theorem \ref{th_Al}. Then there exists $\eps_0>0$ such that, for any $\eps\in(0,\eps_0)$, there exists a smooth solution $u_\eps$ to \eqref{allen-cahn} in $\R^{N}$ such that 
\begin{enumerate}
\item $u_\eps$ is $O(m)\times O(n)$-invariant;

\medskip
\item the nodal set of $u_{\eps}$ has exactly two connected components which are normal graphs over $\Sigma_\eps$, i.e.,
$$
u_{\eps}(\x)= w(z- \h_1(\eps \pe)) - w(z-\h_2(\eps \pe)) -1+o_\eps(1)
$$
where $\h_1,\h_2:\Sigma\to \R$ are $O(m)\times O(n)$-invariant functions.
\medskip
\item $m(u_{\eps})=+\infty$ and for any $\eps >0$ small and any $R> 2\eps^{-1}$,
\begin{equation}\label{est-en}
\int_{B_R}\frac{1}{2} |\nabla u_\eps|^2+\frac{1}{4}(1-u_\eps^2)^2\le c R^N.
\end{equation}
\end{enumerate}
\end{theorem}

In low dimensions \eqref{est-en} suggests stability of the solution according to a conjecture raised in \cite{DKW} (see also \cite{KWang}).  However, Theorem \ref{Theorem 1} provides a counterexample to such conjecture.


On the other hand, the normal graphs $\h_1,\h_2$ are approximate solutions to the Jacobi-Toda system
\begin{equation}\label{JacTodaSystIntro}
\begin{aligned}
\eps^2 \big(\Delta_{\Sigma} \h_1 + |A_{\Sigma}|^2 \h_1\big) - a_0 e^{-\sqrt{2}(\h_2 -\h_1)}&=0\\
\eps^2 \big(\Delta_{\Sigma} \h_2 + |A_{\Sigma}|^2 \h_2\big) + a_0 e^{-\sqrt{2}(\h_2 -\h_1)}&=0
\end{aligned}
\quad \hbox{in} \quad \Sigma,
\end{equation}
where $a_{0}>0$ constant depending only on $w$ (see \eqref{heteroclinic}). By setting ${\tt v}_{1}=\h_1+\h_2$, ${\tt v}_{2}=\h_1-\h_2$, system \eqref{JacTodaSystIntro} decouples into the system
\begin{equation}\label{JacTodaSystIntro2}
\begin{aligned}
\Delta_{\Sigma} {\tt v}_{1} + |A_{\Sigma}|^2 {\tt v}_{1}&=0\\
\eps^2 \big(\Delta_{\Sigma} {\tt v}_{2}+ |A_{\Sigma}|^2 {\tt v}_{2} \big) - 2a_0 e^{-\sqrt{2}{\tt v}_{2}}&=0
\end{aligned}
\quad \hbox{in} \quad \Sigma.
\end{equation}

The first equation in (\ref{JacTodaSystIntro2}) can be easily solved by choosing ${\tt v}_1=0$, i.e., $\h_1=-\h_2$. Concerned with the second one, we present the following result. 

\begin{theorem}[\cite{AKR}]
Under the assumptions of Theorem \ref{Theorem 1} on $\Sigma$, for any $\eps>0$ sufficiently small, the equation
\begin{equation}
\label{eq_Liouville}
\eps^2 \big(\Delta_\Sigma {\tt v}+|A_\Sigma|^2 {\tt v}\big)=  2a_{\star}e^{-\sqrt{2}{\tt v}} \quad \hbox{in }\Sigma
\end{equation}
has a smooth $O(m)\times O(n)$-invariant solution ${\tt v}={\tt v}(\pe)$, with
\begin{equation}\label{asymptw0infty}
{\tt v}(\pe)={\frac{1}{\sqrt{2}}
\log\bigg(\frac{2\sqrt{2}a_\star}{\eps^2}\bigg)-}{\frac{1}{\sqrt{2}}
\log\bigg(|A_{\Sigma}|^{2}(\pe)\bigg)}
{- \frac{1}{\sqrt{2}}\log \bigg(\log \bigg(\frac{2\sqrt{2}a_\star}{\eps^2|A_{\Sigma}|^2(\pe)}\bigg)\bigg)}+ {o(1)}
\end{equation}
for any $\pe\in \Sigma$ and identity \eqref{asymptw0infty} can be differentiated in $\pe$. 
\end{theorem}

We remark that $
o(1)$ is meant uniformly as $\eps \to 0$ and $|\pe|\to \infty$. On the other hand the proof of Theorem \ref{Theorem 1} is mainly based on an infinite dimensional Lyapunov-Schmidt reduction method. Solving the projected equation involves an invertibility theory for the linearized operator of \eqref{eq_Liouville} around the solution \eqref{asymptw0infty}, namely we need to solve equation
\begin{equation*}
\label{eq_lin_JT_Sigma}
\eps^2\big(\Delta_{\Sigma}{\tt v}_1 + |A_{\Sigma}|^2{\tt v}_1\big)+2\sqrt{2}a_\star e^{-\sqrt{2}{\tt v}_0}{\tt v}_1= {\tt f} \quad \hbox{in}\quad\Sigma,
\end{equation*}
where ${\tt f}:\Sigma\to\R$ is $O(m)\times O(n)$-invariant. In this regard we refer the reader to Proposition 3.1 in \cite{AKR} and to our recent works \cite{AR2,AR3} for more details.

Theorem \ref{Theorem 1} was recently generalized in \cite{AR2} to the case where the nodal set has a finite number $k\ge 2$ of connected components. This generalization is rather involved and requires to solve a Jacobi-Toda type system for multiple interacting normal graphs over $\Sigma$.

\begin{theorem}[\cite{AR2}]
\label{Theorem 1.1} Let $m,\,n\ge 2$, $n+m=N+1\ge 8$ and let $\Sigma$ be one of the hypersurfaces constructed in Theorem \ref{th_Al}. Given any $k\in \mathbb{N}$, there exists $\eps_0>0$ such that, for any $\eps\in(0,\eps_0)$, there exists a smooth solution $u_\eps$ to \eqref{allen-cahn} in $\R^{N+1}$ such that 
\begin{enumerate}
\item $u_\eps$ is $O(m)\times O(n)$-invariant;

\medskip
\item nodal set of $u_{\eps}$ has exactly $k$ connected components which are $O(m)\times O(n)$-invariant normal graphs over $\Sigma_{\eps}$

\medskip
\item $m(u_{\eps})=+\infty$ and for any $\eps >0$ small and estimate \eqref{est-en} holds true.
\end{enumerate}
\end{theorem}

{\bf Acknowledgements} O.
Agudelo was supported by the Grant 22-18261S of the Grant Agency of the Czech Republic. M. Rizzi was
partially supported by the German Research Foundation (DFG), project number 62202684.


\begin{thebibliography}{100}

\bibitem{ADW2015} O. Agudelo, M. del Pino, J. Wei, Solutions with multiple catenoidal ends to the Allen-Cahn equation in R3, J. Math.
Pures Appl. (9) 103 (1) (2015) 142-218 

\bibitem{ADW2016} O. Agudelo, M. del Pino, J. Wei, Higher-dimensional catenoid, Liouville equation, and Allen-Cahn equation, Int. Math. Res. Not.(23) (2016) 7051-7102.

\bibitem{AKR} O. Agudelo, M. Kowalczyk, M. Rizzi, Doubling construction for $O(m)\times O(n)$-invariant solutions to the Allen-Cahn equation. (English summary) Nonlinear Anal. 216 (2022), Paper No. 112705, 53 pp.

\bibitem{AR2} O. Agudelo, M. Rizzi, $k$-ended $O(m)\times O(n)$ invariant solutions to the Allen-Cahn equation with infinite Morse index. J. Funct. Anal. 283 (2022), no. 5, Paper No. 109561, 43 pp.

\bibitem{AR3} O. Agudelo, M. Rizzi, The Jacobi operator of some spcial minimal hypersurfaces, arXiv:2408.08728.

\bibitem{ABPRS} H. Alencar, A. Barros, O. Palmas, G.  Reyes, W. Santos, $O(m)\times O(n)$-invariant minimal hypersurfaces in $\R^{m+n}$. \emph{Ann. Global Anal. Geom.} {27} (2005), no. 2, 179--199.


\bibitem{Allen-Cahn} S. Allen, J.W. Cahn, A microscopic theory for antiphase boundary motion and its application to antiphase domain coarsening, Acta Metall. 27 (1979) 1084-1095.


\bibitem{AC} L. Ambrosio, X. Cabr\'{e}, Entire solutions of semilinear elliptic equations in $\R^3$ and a conjecture of De Giorgi. J. Amer. Math. Soc. 13 (2000), no. 4, 725-739.

\bibitem{Cabre2011} X. Cabr\'{e}
Uniqueness and stability of saddle-shaped solutions to the Allen-Cahn equation
J. Math. Pures Appl. (9), 98 (3) (2012), pp. 239-256, 10.1016/j.matpur.2012.02.006

\bibitem{CT2009} X. Cabr\'{e}, J. Terra, Saddle-shaped solutions of bistable diffusion equations in all of $R^{2m}$, J. Eur. Math. Soc. 11(4) (2009) 819-843.

\bibitem{CT2010} X. Cabr\'{e}, J. Terra, Qualitative properties of saddle-shaped solutions to bistable diffusion equations, Commun. Partial Differ. Equ. 35(11) (2010) 1923-1957.


\bibitem{DFP} Dang, H., Fife, P.C.Peletier, L.A. Saddle solutions of the bistable diffusion equation. Z. angew. Math. Phys. 43, 984-998 (1992). https://doi.org/10.1007/BF00916424

\bibitem{DKW-degiorgi} M. del Pino, Manuel, M. Kowalczyk, J. Wei, On De Giorgi's conjecture in dimension $N\ge 9$. (English summary)
Ann. of Math. (2) 174 (2011), no. 3, 1485-1569.

\bibitem{DKW} M. del Pino, M. Kowalczyk, J. Wei,  Entire solutions of the Allen-Cahn equation and complete embedded minimal surfaces of finite total curvature in R3. J. Differential Geom. 93 (2013), no. 1, 67--131.

\bibitem{GG} N. Ghoussoub, C. Gui, On a conjecture of De Giorgi and some related problems. Math. Ann. 311 (1998), no. 3, 481-491.

\bibitem{GLW} C. Gui, K. Wang, J. Wei, Axially symmetric solutions of the Allen-Cahn equation with finite Morse index (English summary), Trans. Am. Math. Soc. 373(5) (2020) 3649-3668.


\bibitem{KS} R. Kohn, Robert, P. Sternberg, Local minimisers and singular perturbations.
Proc. Roy. Soc. Edinburgh Sect. A 111 (1989), no. 1-2, 69-84.

\bibitem{M} L. Mazet, Minimal hypersurfaces asymptotic to Simons cones. J. Inst. Math. Jussieu 16 (2017), no. 1, 39--58.

\bibitem{Mo} L. Modica, The gradient theory of phase transitions and the minimal interface criterion. Arch. Rational Mech. Anal. 98 (1987), no. 2, 123-142.

\bibitem{Nunes} I. P. Nunes, Rigidity of Area-minimising hyperbolic surfaces in three manifolds, doctoral thesis, Rio de Janeiro, Instituto de Matematica pura e aplicada.

\bibitem{PR} F. Pacard, M. Ritor\'{e}, From constant mean curvature hypersurfaces to the gradient theory of phase transitions. (English summary)
J. Differential Geom. 64 (2003), no. 3, 359-423.

\bibitem{Sa} O. Savin, Regularity of flat level sets in phase transitions. Ann. of Math. (2) 169 (2009), no. 1, 41-78.

\bibitem{KWang} K. Wang, Some remarks on the structure of finite Morse index solutions to the Allen-Cahn equation in R2, Nonlinear Differ. Equ. Appl. 24(5) (2017) 58, 17 pp.

\end{thebibliography}
\end{document}